\documentclass{article}

\usepackage{amsmath,amssymb,latexsym,epsfig,amsfonts}

\newtheorem{theorem}{Theorem}[section]
\newtheorem{proposition}[theorem]{Proposition}
\newtheorem{corollary}[theorem]{Corollary}

\newtheorem{lemma}[theorem]{Lemma}

\newcommand{\beha}{\begin{enumerate}}
\newcommand{\behe}{\end{enumerate}}

\newcommand{\bR}{{\bf R}}
\newcommand{\bC}{{\mathbb C}}

\newcommand{\bN}{{\bf N}}

\newcommand{\RE}[1]{(\ref{#1})}
\newcommand{\DD}{{\rm d}}
\newcommand{\eproof}{$\quad\Box$\\[0.3cm]}

\begin{document}
\title{On measures of maximal and full dimension for  polynomial automorphisms of $\bC^2$}
\renewcommand{\baselinestretch}{1.0}
\author{
Christian Wolf \footnote{The author was supported by a research
fellowship of the Deutsche Forschungsgemeinschaft (DFG)}
\\
 Department of Mathematics, Indiana University\\
 Bloomington, IN 47405, USA\\
 Email: cwolf@indiana.edu
}
 \maketitle
\begin{abstract}
For a hyperbolic polynomial automorphism of $\bC^2$, we show the
existence of a measure of  maximal dimension, and identify the
conditions under which a measure of full dimension
exists.\\[0.6cm] AMS 2000 Subject Classification: 37FXX, 37C45,
32HXX
\end{abstract}
\section{Introduction}
Let $g$ be a hyperbolic polynomial automorphism of $\bC^2$. An
important approach for understanding the dynamics of $g$ is to
study its invariant Borel probability measures. One key idea is to
study the Hausdorff dimensions of these measures. For an invariant
Borel probability measure $\nu$, we define the Hausdorff dimension
of $\nu$ by
\begin{equation}\label{eqDD}
\dim_H (\nu)=\inf\{\dim_H A:  \nu(A)=1\}.
\end{equation}
We define the quantity $\DD(g)$ by
\begin{equation}\label{eqqDDD}
\DD(g)=\sup\{\dim_H(\nu)\},
\end{equation}
where the supremum is taken over all ergodic invariant  Borel
probability measures with positive entropy.  This quantity was
introduced by Denker and Urbanski \cite{DU} in the context of
rational maps on the Riemann sphere. They called it the dynamical
dimension of the map.

 It is easy to see that the support of each
measure considered in \RE{eqqDDD}  is contained in the Julia set
$J$ (see Section 2 for the definition). We denote by $M(J,g|_J)$
the set of all ergodic invariant  Borel probability measures
supported on $J$.
 If $\nu\in M(J,g|_J)$, then Young's formula \cite{Y} implies that  $h_\nu(g)>0$ if and only if
$\dim_H(\nu)> 0$. Thus,  we could have also included measures with
zero entropy while taking the supremum in \RE{eqqDDD}.

If a measure $\nu\in M(J,g|_J)$ attains the supremum in
\RE{eqqDDD}, that is,
\begin{equation}
\dim_H(\nu)=\DD(g),
\end{equation}
we say that $\nu$ is a measure of maximal dimension for $g$.

 McCluskey and Manning  \cite{MM} mentioned a heuristic argument for the
existence of a measure of maximal dimension in the case of Axiom A
surface diffeomorphisms. However, it is until today not known
whether this argument can be extended to a rigorous proof (see the
remarks after Theorem \ref{mesmaxdim} for more details).

In this paper we study the existence of a measure of maximal
dimension for hyperbolic polynomial automorphisms of $\bC^2$. Our
main result is the following.
\begin{theorem}\label{thmainresult} Let $g$ be a hyperbolic polynomial automorphism of $\bC^2$. Then
there exists a measure of maximal dimension for $g$. The set of
measures of maximal dimension  is finite.
\end{theorem}
The proof of this theorem uses the theory of thermodynamic
formalism. The key idea is to extract a one-parameter family of
potentials and to consider the corresponding family of equilibrium
measures. We show that a measure of maximal dimension necessarily
belongs to  this family of equilibrium measures. Furthermore, if a
measure $\nu$ maximizes Hausdorff dimension among these
equilibrium measures, then $\nu$ is a measure of maximal
dimension.

The Hausdorff dimension of the Julia set is also an interesting
dimension-theoretical feature of $g$. Recently, there has been
made substantial progress on this subject (see \cite{FO},
\cite{FS}, \cite{VW}, \cite{W1}, \cite{W2}). We say that $\nu\in
M(J,g|_J)$ is a measure of full dimension if
\begin{equation}
\dim_H(\nu)=\dim_H J.
\end{equation}
 Friedland and Ochs \cite{FO} studied
 the existence of a measure of full dimension.
 They proved existence  for volume preserving maps.
 We provide an alternative proof for this result in Theorem
\ref{thFO}. Moreover, we obtain in Corollary \ref{thmefulun} that
there exists at most one measure of full dimension. In the next
theorem we consider non-volume preserving maps.
\begin{theorem}\label{th2}
Let $g$ be a non-volume preserving hyperbolic polynomial
automorphism of $\bC^2$. Assume that $\nu$ is a measure of full
dimension for $g$. Then $\nu$ is the measure of maximal entropy of
$g$ .
\end{theorem}
Theorem \ref{th2} indicates that for non-volume preserving maps
the existence of a measure of full dimension  seems to be a very
rare phenomenon. Indeed we do not have an example of a non-volume
preserving hyperbolic polynomial automorphism of $\bC^2$ admitting
a measure of full dimension.

In Section 4 we observe that there exists a dense open subset of
hyperbolic parameter space for which no measure of full dimension
exists. This implies that
\begin{equation}
\DD(g)<\dim_H J
\end{equation}
holds for  these parameters.

In the last part of this paper we analyze the dependence of
$\DD(g)$ on the parameter of the mapping. More precisely, we prove
the following result.
\begin{theorem} Let $\lambda\mapsto g_\lambda$ be a holomorphic family of hyperbolic polynomial
automorphisms of $\bC^2$ of fixed degree. Then $\lambda\mapsto
\DD(g_\lambda)$ is continuous and plurisubharmonic.\end{theorem}
This paper is organized as follows. In Section 2 we present the
basic definitions and notations. In Section 3 we introduce
elements from dimension theory for hyperbolic polynomial
automorphisms of $\bC^2$ and provide the tools  for the analysis
of the existence of measures of maximal and full dimension.
Section 4 is devoted to the analysis of the existence of a measure
of full dimension. The existence of a measure of maximal dimension
is proved in Section 5. Finally, we study in Section 6 the
dependence of $\DD(g)$ on parameters.\\[0.3cm]

It would be interesting to understand whether, or at least under
which conditions, an uniqueness result for the measure of maximal
dimension holds. A partial answer to this question is given in
Corollary \ref{thmefulun} where an uniqueness result is shown in
the case when a measure of full dimension exists.

\section{Notation and Preliminaries}
Let $g$ be a polynomial automorphism of $\bC^2$. We can associate
with $g$ a dynamical degree $d$ which is a conjugacy invariant
(see \cite{FM}, \cite{BS2} for details). We are interested in
nontrivial dynamics which  occurs if and only if $d>1$. Friedland
and Milnor  showed in \cite{FM} that every polynomial automorphism
of $\bC^2$ with nontrivial dynamics  is conjugate to a mapping of
the form $g=g_1\circ...\circ g_m$, where each $g_i$ is a
generalized H\'enon  mapping. This means that $g_i$ has the form
\begin{equation}\label{defhen}
g_i(z,w)= (w,P_i(w)+a_i z),
\end{equation}
where $P_i$ is a complex polynomial of degree  $d_i\geq 2$ and
$a_i$ is a non-zero complex number. The dynamical degree $d$ of
$g$ is equal to $d_1\cdot...\cdot d_m$ and therefore coincides
with the polynomial degree of $g$. In this paper we assume that
$g$ is a finite composition of generalized H\'enon mappings.
 Since dynamical properties are invariant under conjugacy,
 our results also hold for a general polynomial automorphism
of $\bC^2$ with nontrivial dynamics.

For $g$ we define $K^\pm$ as the set of points in $\bC^2$ with
bounded forward/backward orbits, $K=K^+\cap K^-$, $J^\pm=
\partial K^\pm$ and $J= J^+\cap J^-$.
 We refer to $J^\pm$ as the
positive/negative Julia set of $g$ and $J$ is the Julia set of
$g$. The sets $K$ and $J$ are compact.

Note that the function $a=\det Dg$ is constant in $\bC^2$.
Therefore we can restrict our considerations to the volume
decreasing case $(|a|<1)$, and to the volume preserving case
$(|a|=1)$, because otherwise we can consider $g^{-1}$. This will
be a standing assumption in this paper. We note that for $g$ we
have $a=a_1\cdot...\cdot a_m$.

As pointed out in the introduction we assume in this paper that
$g$ is a hyperbolic mapping. This means that  there exists a
continuous invariant splitting $T_J \bC^2= E^u\oplus E^s$ such
that $Dg|_{E^u}$ is uniformly expanding and $Dg|_{E^s}$ is
uniformly contracting. Hyperbolicity implies that we can associate
with each point $p\in J$ its local unstable/stable manifold
$W^{u/s}_\epsilon(p)$. Moreover $g$ is an Axiom A diffeomorphism
(see \cite{BS1} for more details).
\section{Elements from dimension theory}
In this section we introduce elements from dimension theory for
hyperbolic polynomial automorphisms of $\bC^2$ and provide the
tools for the analysis of measures of maximal and full dimension.

We start   by introducing Lyapunov exponents. Let $\nu\in
M(J,g|_J)$. By the multiplicative ergodic Theorem of Oseledec,
there are Lyapunov exponents $\lambda^{(1)}_\nu\leq
\lambda^{(2)}_\nu$ with respect to $\nu$ (see e.g. \cite{KH}). The
Julia set $J$ is a hyperbolic set of index 1 (see \cite{BS1}).
This implies
\begin{equation}
\lambda^{(1)}_\nu<0<\lambda^{(2)}_\nu.
\end{equation}
In particular, $\nu$ is a hyperbolic measure. We define the
quantity
\begin{equation}\label{poslya}
\Lambda(\nu)=\lim_{n\to\infty} \frac{1}{n} \int \log||Dg^n||d\nu.
\end{equation}
Similar as it was done for the measure of maximal entropy in
\cite{BS3}, the positive Lyapunov exponent $\lambda^{(2)}_\nu$
coincides with $\Lambda(\nu)$. Since $g$ has constant jacobian
determinant $a$, the negative Lyapunov exponent
$\lambda^{(1)}_\nu$ is given by $-\Lambda(\nu)+\log|a|$.
 In \cite{W1} we derived the formula
\begin{equation}\label{lamw3}
\Lambda(\nu)=\int \log ||Dg|_{E^u}||d\nu.
\end{equation}
By Young's formula \cite{Y}, we have for all $\nu\in M(J,g|_J)$
that
\begin{equation}\label{yo}
\dim_H(\nu)=\frac{h_\nu(g)}{\Lambda(\nu)}+\frac{h_\nu(g)}{\Lambda(\nu)-\log
|a|}.
\end{equation}
Here $h_\nu(g)$ denotes the measure theoretic entropy of $g$ with
respect to $\nu$.

Next we introduce topological pressure. Let $C(J,\bR)$ denote the
Banach space of all continuous functions from $J$ to $\bR$. The
topological pressure of $g|_J$, denoted by $P=P(g|_J,.)$, is a
mapping from $C(J,\bR)$ to $\bR$ (see \cite{W} for the
definition). We have $P(g|_J,0)= h_{top}(g)=\log d$, where
$h_{top}(g)$ denotes the topological entropy of $g$ and $d$ is the
polynomial degree of $g$. The variational principle provides the
formula
\begin{equation}\label{eqvarpri}
P(g|_J,\varphi)= \sup_{\nu\in M(J,g|_J)} \left(h_\nu(g)+\int_J
\varphi d\nu\right).
\end{equation}
  If a measure $\nu_\varphi\in M(J,g|_J)$
achieves the supremum in equation \RE{eqvarpri}, that is,
\begin{equation}\label{eqequil}
P(g|_J,\varphi)= h_{\nu_\varphi} (g)+\int_J\varphi d\nu_\varphi,
\end{equation}
it is called  equilibrium measure of the potential $\varphi$.

The topological pressure has the following properties (see
\cite{R}). \beha
\item  The topological pressure  is a convex function.
\item If $\varphi$ is a strictly negative function, then the
function $t\mapsto P(g|_J,t\varphi)$ is strictly decreasing.
\item The topological pressure is a real analytic function on the subspace of H\"{o}lder
continuous functions, that is, when $\alpha >0$ is fixed, then
 $P(g|_J,.)|_{C^\alpha(J,\bR)}$ is  a real analytic function. Note that
  $C^\alpha$ can not be replaced by $C^0$.
 \item If $\alpha>0$ and $\varphi\in C^\alpha(J,\bR)$, then there exists a uniquely defined
equilibrium measure $\nu_\varphi\in M(J,g|_J)$ of the potential
$\varphi$.
  Furthermore we have for all $\varphi, \psi\in C^\alpha(J,\bR)$
 \begin{equation}\label{eqdifpre}
  \frac{d}{dt}\bigg|_{t=0} P(g|_J,\varphi + t\psi )= \int_\Lambda \psi
  d\nu_\varphi.
 \end{equation}\behe
We now introduce  potentials which are related to  Lyapunov
exponents. We define
\begin{equation}
\phi^{u/s}:J\to \bR,\qquad p\mapsto \log ||Dg(p)|_{E^{u/s}_p}||
\end{equation}
and  the unstable/stable pressure function
\begin{equation}
P^{u/s}:\bR\to\bR,\qquad t\mapsto P(g|_J,\mp t\phi^{u/s}).
\end{equation}
Since $\phi^{u/s}$ is H\"{o}lder continuous (see \cite{B}), we may
follow from property iii) of the topological pressure that
$P^{u/s}$ is real analytic. Property iv) of the topological
pressure implies that there exists a uniquely defined  equilibrium
measure $\nu_{\mp t\phi^{u/s}}\in M(J,g|_J)$ of the potential $\mp
t\phi^{u/s}$.

We will need the following result about the relation between the
unstable and stable pressure function.
\begin{proposition}[\cite{W2}]\label{thtp}
$P^u(t)=P^s(t)-t\log |a|.$
\end{proposition}
\begin{lemma}\label{propu=s}
$\nu_{-t\phi^{u}}=\nu_{t\phi^s}$
\end{lemma}
{\it Proof. } Let $t\geq 0$. Then
\begin{eqnarray}
 P^s(t)&=& P^u(t)+t\log|a|\nonumber\\ &=& h_{\nu_{-t\phi^u}}(g)
+t\left(-\int \log ||Dg|_{E^u}||d\nu_{-t\phi^u}
+\log|a|\right)\nonumber\\ &=& h_{\nu_{-t\phi^u}}(g)
+t\left(\lim_{n\to\infty}\frac{1}{n}\int-\log||Dg^{n}|_{E^u}||
+\log|a^n|d\nu_{-t\phi^u}\right)\nonumber\\ &=&
h_{\nu_{-t\phi^u}}(g) +t\left(\lim_{n\to\infty}\frac{1}{n}\int
\log||Dg^{n}|_{E^s}|| d\nu_{-t\phi^u}\right)\nonumber\\
 &=&
h_{\nu_{-t\phi^u}}(g) +t\int \log
||Dg|_{E^s}||d\nu_{-t\phi^u}\nonumber\\ &=& h_{\nu_{-t\phi^u}}(g)
+t\int \phi^s d\nu_{-t\phi^u}.
\end{eqnarray}
The result follows from the uniqueness of the equilibrium measure
of the potential $t\phi^s$.\eproof We will use in the remainder of
this paper the notation $\nu_t=\nu_{\mp t\phi^{u/s}}$. This
notation is justified by Lemma \ref{propu=s}. We also write
$\Lambda(t)=\Lambda(\nu_t)$ and $h(t)=h_{\nu_t}(g)$, and consider
$\Lambda$ and $h$ as real-valued functions of $t$. Equations
\RE{lamw3}, \RE{eqequil} imply
\begin{equation}\label{eqphl}
P^u(t)=h(t)-t\Lambda(t).
\end{equation}
Therefore, Proposition \ref{thtp} implies
\begin{equation}\label{eqphs}
P^s(t)=h(t)-t(\Lambda(t)-\log|a|).
\end{equation}

\begin{proposition}\label{proisnul}
$\Lambda$ and $h$ are real analytic. Furthermore
\begin{equation}
\frac{d\Lambda}{dt}\leq 0.
\end{equation}
 If $\Lambda$ is not constant, then every zero of the  derivative of $\Lambda$ is
 isolated.
\end{proposition}
{\it Proof. } Let $t_0\geq 0$. We define potentials
$\varphi=-t_0\phi^u, \psi=-\phi^u$. Therefore application of
equations \RE{lamw3} and    \RE{eqdifpre} imply
\begin{equation}\label{eq1}
\frac{dP^u}{dt} \left(t_0\right) =
-\Lambda(\nu_{t_0})=-\Lambda(t_0).
\end{equation}
Since $P^u$ is  real analytic, we obtain that $\Lambda$ is  also
real analytic. We conclude from \RE{eqphl} that $h$ is also  real
analytic. The convexity of $P^u$ implies
\begin{equation}
\frac{d^2P^u}{dt^2} \geq 0,
\end{equation}
hence
\begin{equation}
\frac{d\Lambda}{dt}  \leq 0.
\end{equation}
Finally, if $\Lambda$ is not constant, then the uniqueness theorem
for real analytic functions, applied to the derivative of
$\Lambda$, implies that all zeros of the derivative of $\Lambda$
are isolated.\eproof
\begin{corollary}\label{corcondec}
$\Lambda$ is either constant or strictly decreasing.
\end{corollary}
{\it Proof. } The result follows immediately from Proposition
\ref{proisnul}.\eproof Of particular interest for the analysis of
the existence of measures of maximal and full dimension are the
Hausdorff dimensions of the measures $\nu_t$. We use the notation
$\Delta(t)=\dim_H(\nu_t)$. Equation \RE{yo} yields
\begin{equation}\label{defdelta}
\Delta(t)=\frac{h(t)}{\Lambda(t)}+\frac{h(t)}{\Lambda(t)-\log
|a|}.
\end{equation}
Thus, $\Delta$ is also a real analytic function. Equations
\RE{eqphl}, \RE{eqphs} and Proposition \ref{thtp} imply
\begin{equation}\label{eqfordim}
\Delta(t)=2t+\frac{P^u(t)}{\Lambda(t)}+\frac{P^u(t)+t\log|a|}{\Lambda(t)-\log|a|}.
\end{equation}
From an elementary calculation we obtain the following formula for
the derivative of $\Delta$.
\begin{equation}\label{eqdiffdim}
\frac{d\Delta}{dt}(t_0)=-\frac{\frac{d\Lambda}{dt}(t_0)\left[P^u(t_0)(\Lambda(t_0)-\log|a|)^2
+ (P^u(t_0)+t_0\log|a|)\Lambda(t_0)^2\right]}{\Lambda(t_0)^2
(\Lambda(t_0)-\log|a|)^2}.
\end{equation}

Finally we consider the Hausdorff dimension of the Julia set $J$.
The following result due to Verjovsky and Wu  provides a formula
for the Hausdorff dimension of the unstable/stable slice in terms
of the zeros of the pressure functions.
\begin{theorem}[\cite{VW}]\label{thVW}
Let $p\in J$. Then $t^{u/s}=\dim_H W^{u/s}_\epsilon(p)\cap J$ does
not depend on $p\in J$. Furthermore $t^{u/s}$ is given by the
unique solution of
\begin{equation}\label{eqBR}
P^{u/s}(t)=0.
\end{equation}
\end{theorem}
Equation \RE{eqBR} is called Bowen-Ruelle formula. We refer to
$t^{u/s}$ as the Hausdorff dimension of the unstable/stable slice.

In \cite{W1} we proved the formula
\begin{equation}\label{djtuts} \dim_H J = t^u+t^s=\sup_{\nu\in M(J,g|_J)}\left(\frac{h_\nu(g)}{\Lambda(\nu)}\right)+
\sup_{\nu\in
M(J,g|_J)}\left(\frac{h_\nu(g)}{\Lambda(\nu)-\log|a|}\right),
\end{equation}
where each of the suprema on the right-hand side of the equation
is uniquely attained by the measures $\nu_{t^u}$ and $\nu_{t^s}$
respectively. Hence
\begin{equation}\label{eqmesdimJ}
\dim_H J = \frac{h(t^u)}{\Lambda(t^u)} +
\frac{h(t^s)}{\Lambda(t^s)-\log|a|}.
\end{equation}
Equation \RE{eqmesdimJ} and the uniqueness of the measures
$\nu_{t^u}, \nu_{t^s}$ in equation \RE{djtuts} imply that, if
there exists a measure of full dimension, then it already
coincides with $\nu_{t^u}$ and $\nu_{t^s}$. Thus, we have the
following result.
\begin{corollary}\label{thmefulun}
 Assume $m$
is a measure of full dimension for $g$, then $m=\nu_{t^u}=
\nu_{t^s} $. In particular, there exists at most one measure of
full dimension.
\end{corollary}

\section{Measures of full dimension}
In this  section we identify the conditions under which a measure
of full dimension exists. More precisely, we show that a measure
of full dimension exists if and only if $g$ is either volume
preserving, or $g$ is volume decreasing and the measure of maximal
entropy is a measure of full dimension.

We start with the volume preserving case.
\begin{theorem}\label{thFO}
Let $g$ be volume preserving. Then $t^u=t^s$, and $\nu_{t^u}$ is a
measure of full dimension for $g$.
\end{theorem}
{\it Proof. } We have $|a|=1$, therefore Proposition \ref{thtp}
implies $P^u=P^s$. Thus, Theorem \ref{thVW} yields $t^u=t^s$.
Therefore, by equations \RE{yo}, \RE{eqmesdimJ}, we conclude that
$\dim_H (\nu_{t^u})=\dim_H J$, which implies that $\nu_{t^u}$ is a
measure of full dimension. \eproof {Remark. } As noted in the
introduction, in the volume preserving case the existence of a
measure of full dimension was already shown by Friedland and Ochs
\cite{FO}. They proved that the existence of a measure of full
dimension follows from the fact that $|\det Dg^n(p)|=1$ holds for
every periodic point $p$ with period $n$. They also observed that
in this case the measure of full dimension is equivalent to the
$t$-dimensional Hausdorff measure, where $t$ is the Hausdorff
dimension of $J$.
\\[0.3cm] We now consider the volume decreasing case. The following theorem is the main result of this
section.
\begin{theorem}\label{thmesmax} Assume $g$ is volume decreasing. Then the following
are equivalent. \beha
\item $g$ admits a  measure of full dimension.
\item The unstable pressure function $P^u$ is affine.
\item The stable pressure function $P^s$ is affine.
\item The measure of maximal entropy is a measure of full
dimension for $g$. \behe
\end{theorem}
{\it Proof.}\\
 $ii) \Leftrightarrow iii)$ follows from Proposition
\ref{thtp}.\\ $i)\Rightarrow ii)$ Let us assume that $g$ admits a
measure of full dimension. Corollary \ref{thmefulun} implies that
$\nu_{t^u}=\nu_{t^s}$ is the measure of full dimension. Since
$t^s<t^u$ (see \cite{W2}), Corollary \ref{corcondec} implies that
$P^u$ has constant derivative in $[t^s,t^u]$. Therefore, since
$P^u$ is real analytic, we may conclude  that $P^u$ is affine.\\
$ii)+iii) \Rightarrow iv)$ Let $\mu$ denote the measure of maximal
entropy. The topological entropy of $g|_J$ is equal to $\log d$
(see \cite{BS3}). Thus $P^u(0)=P^s(0)=\log d$. Equation \ref{eq1}
and Proposition \ref{thtp} imply
\begin{eqnarray}
\frac{dP^u}{dt} (0)&=&-\Lambda(\mu)
\\
\frac{dP^s}{dt} (0)&=&-\Lambda(\mu)+\log|a|.
\end{eqnarray}
Since $P^u$ and $P^s$ are affine,  Theorem \ref{thVW} and equation
\RE{djtuts} imply
\begin{equation}\label{eq3}
\dim_H J = \frac{\log d}{\Lambda(\mu)} + \frac{\log
d}{\Lambda(\mu)-\log|a|}.
\end{equation}
But by Young's formula \RE{yo}, the right-hand side of equality
\RE{eq3} is equal to $\dim_H(\mu)$. Thus, $\mu$ is a measure of
full dimension.\\ $iv)\Rightarrow i)$ trivial. \eproof
\begin{corollary}\label{corconcan}Assume $g$ is volume decreasing. If there exists
a
measure of full dimension for $g$, then $J$ is either a Cantor set
or connected.
\end{corollary}
{\it Proof.} Let $\mu$ denote the measure of maximal entropy. It
is shown in \cite{BS3} that $\Lambda(\mu)\geq \log d$. The Julia
set $J$ is connected, if and only if $\Lambda(\mu)=\log d$  (see
\cite{BS6}). Therefore, if $J$ is not connected, then
$\Lambda(\mu)>\log d$. Thus, equation \RE{djtuts} implies
$t^u,t^s<1$. Since $J$ has a local product structure (see
\cite{BS1}), we may conclude that $J$ is a Cantor set. \eproof Let
us assume that $g$ is volume decreasing and let $\mu$ denote the
measure of maximal entropy of $g$. We assume that $g$ admits a
measure of full dimension; thus,  by Theorem \ref{thmesmax}, $\mu$
is the measure  of full dimension for $g$. Let ${\cal S}$ denote
the set of all saddle points of $g$. For $p\in{\cal S}$ with
period $n$ we denote by $\lambda^{u/s}(p)$ the eigenvalues of
$Dg^n(p)$, where $|\lambda^s(p)|<1<|\lambda^u(p)|$. Then, we may
conclude by Theorem \ref{thmesmax} and  Proposition 4.5 of
\cite{B} that
\begin{equation}\label{equaleig}
\log |\lambda^u(p)|= n\Lambda(\mu)
\end{equation}
holds for all $p\in {\cal S}$.
 Using a perturbation
argument, it is not to hard to see that we can find arbitrary
close to $g$ a hyperbolic polynomial automorphism of $\bC^2$ for
which \RE{equaleig} does not hold, and which therefore does not
admit a measure of full dimension. Here we mean close with respect
to the topology on hyperbolic parameter space induced by the
parameter of the mapping (see \cite{W1} for details). On the other
hand, it is obvious that the set of parameters having no measure
of full dimension is an open subset of hyperbolic parameter space.
Thus, there exists a dense open subset of parameters admitting no
measure of full dimension. We leave the details to the
reader.\\[0.3cm] { Remarks. } It is a well-known fact that for
hyperbolic and parabolic rational maps on the Riemann sphere the
Hausdorff dimension of the Julia set can be represented in terms
of the Bowen-Ruelle formula (see for instance \cite{U} and the
references therein). Therefore, for these maps the existence of a
measure of full dimension follows  as a consequence of the
variational principle. For results concerning the existence of
measures of full dimension for other maps see \cite{GP} and the
references therein.
\section{Measures of maximal dimension}
In this section we consider the case when $g$ does not admit a
measure of full dimension. Thus, by the results of the previous
section, we may assume that  $g$ is volume decreasing and
$P^{u/s}$ is not affine. This will be a standing assumption in
this section. The following theorem is the main  result of this
paper.
\begin{theorem}\label{mesmaxdim}
There exists a measure of maximal dimension for $g$. If $m$ is a
measure  of maximal dimension for $g$, then there exists $t^s < t
< t^u$ such that $m$ is the equilibrium measure of the potential
$-t\phi^u$, that is $m=\nu_t$.
\end{theorem}
{\it Proof. } Since $g$ is volume decreasing, we have $t^s<t^u$
(see \cite{W2}). \\ {\it Assertion 1. } There exists $\epsilon>0$
such that $\Delta$ is strictly increasing on $[0,t^s+\epsilon)$
and strictly decreasing on $(t^u-\epsilon,\infty)$.

{\it Proof of Assertion 1. } Theorem \ref{thVW} and the fact that
$P^{u/s}$ is a strictly decreasing function [property ii) of the
topological pressure] imply that  $P^s(t)>0$ for all $t\in
[0,t^s)$. Analogously we have $P^u(t)>0$ for all $t\in [0,t^u)$.
We conclude from Proposition \ref{proisnul}, equation
\RE{eqdiffdim} and an elementary continuity argument that there
exists $\epsilon>0$ such that
\begin{equation}
\frac{d\Delta}{dt}\geq 0
\end{equation}
in $[0,t^s+\epsilon)$, and all zeros of the derivative of $\Delta$
in $[0,t^s+\epsilon)$ are isolated. Therefore, $\Delta$ is
strictly increasing in $[0,t^s+\epsilon)$. A similar argumentation
shows that there exists $\epsilon>0$ such that $\Delta$ is
strictly decreasing in $(t^u-\epsilon,\infty)$.\\ Assertion 1
implies that there exists $t^*\in [t^s+\epsilon,t^u-\epsilon]$
such that
\begin{equation} \dim_H (\nu_{t^*})= \sup_{t\geq 0} \Delta(t).
\end{equation}
{\it Assertion 2.} The measure $\nu_{t^*}$ is a measure of maximal
dimension.

 {\it Proof of Assertion 2. } Let $(m_k)_{k\in
\bN}$ be a sequence in $M(J,g|_J)$ such that
\begin{equation}
\lim_{k\to\infty} \dim_H(m_k) = \DD(g).
\end{equation}
By Assertion 1, we may assume without loss of generality that
$\dim_H (\nu_0)=\Delta(0) < \dim_H (m_k)$ for all $k\in \bN$.
Since $\nu_0$ is the unique measure of maximal entropy, we may
conclude by Young's formula \RE{yo} that
\begin{equation}\label{eqlanuk}
\Lambda(\nu_0)>\Lambda(m_k)
\end{equation}
for all $k\in \bN$. Again by Assertion 1, we may assume without
loss of generality that $\dim_H(\nu_{t^u})=\Delta(t^u) < \dim_H
(m_k)$ for all $k\in \bN$. Equation \RE{djtuts} implies
\begin{equation}\label{eqtunk}
\frac{h_{m_k}(g)}{\Lambda(m_k)}< \frac{h(t^u)}{\Lambda(t^u)}
\end{equation}
for all $k\in \bN$. Therefore, Young's formula \RE{yo} implies
\begin{equation}\label{eqasd}
\frac{h_{m_k}(g)}{\Lambda(m_k)-\log|a|}>
\frac{h(t^u)}{\Lambda(t^u)-\log|a|}
\end{equation}
for all $k\in \bN$. Equations \RE{eqtunk}, \RE{eqasd} imply
$h_{m_k}(g)>h(t^u)$, and therefore again by equation \RE{eqtunk}
we obtain
\begin{equation}\label{eqkktu}
\Lambda(m_k)>\Lambda(t^u)
\end{equation}
for all $k\in \bN$. Since $\Lambda$ is  continuous, equations
\RE{eqlanuk}, \RE{eqkktu} imply that  for all $k\in\bN$ there
exists $t_k\in (0,t^u)$ such that
\begin{equation}
\Lambda(m_k)=\Lambda(t_k).
\end{equation}
Thus, the variational principle \RE{eqvarpri} implies
\begin{equation}
 h_{m_k}(g)\leq h(t_k),
\end{equation}
hence
\begin{equation}
\dim_H(m_k)\leq \Delta(t_k)
\end{equation}
for all $k\in \bN$. This implies
\begin{equation}
\dim_H(m_k)\leq \dim_H(\nu_{t^*})
\end{equation}
for all $k\in \bN$. We conclude that $\nu_{t^*}$ is a measure of
maximal dimension.\\ To complete the proof of the theorem we have
to show the following.\\ {\it Assertion 3.} For every measure $m$
of maximal dimension there exists $t^s < t < t^u$ such that $m$ is
the equilibrium measure of the potential $-t\phi^u$.

{\it Proof of Assertion 3. } Let $m$ be a measure of maximal
dimension. We apply  to $m$ (instead of $m_k$) the same
argumentation as in the proof of Assertion 2. This implies that
there exists  $t\in (0,t^u)$ such that $\Lambda(m)=\Lambda(t)$.
Since $\dim_H(m)\geq \Delta(t)$, we may follow  by equation
\RE{yo}  that $h_{m}(g)\geq h(t)$. On the other hand, since
$\nu_t$ is the equilibrium measure of the potential $-t\phi^u$, we
may conclude by \RE{eqvarpri}, \RE{eqequil} that $h_{m}(g)\leq
h(t)$. Hence $h_{m}(g)= h(t)$. Therefore, the uniqueness of the
equilibrium measure of the potential $-t\phi^u$ implies $m=\nu_t$.
Finally, Assertion 1 implies that $t\in (t^s,t^u)$. This completes
the proof.\eproof { Remarks. } The following heuristic argument
was mentioned by McCluskey and Manning  \cite{MM} to state the
existence of a measure of maximal dimension in the case of $C^2$
axiom A diffeomorphisms of real surfaces. Since the entropy map is
upper semi-continuous it can be shown that the map $\nu\mapsto
dim_H(\nu)$, defined on the set of all ergodic invariant measures,
is also upper semi-continuous. It is now suggested in \cite{MM}
that this implies the existence of a measure of maximal dimension.
To make this argument rigorous we need to show that there exists a
sequence of ergodic invariant measures $m_k$ with
$\dim_H(m_k)\to\DD(g)$ having an ergodic weak* limit. Whether this
holds is not clear since the set of all ergodic invariant measures
is not closed with respect to the weak* topology. The latter
follows for instance from Proposition 21.9 in \cite{DGS}.
\\[0.3cm] We obtain a  formula for $\DD(g)$.
\begin{corollary}
Let $t>0$ such that $\nu_t$ is a measure of maximal dimension.
Then
\begin{equation}
{\DD}(g)= 2t + \frac{P^u(t) \log|a|}{\Lambda(t)^2}.
\end{equation}
\end{corollary}
{\it Proof. }By Proposition \ref{proisnul}, equation
\RE{eqdiffdim} and Theorem \ref{mesmaxdim}, a necessary condition
for $\nu_t$ being a measure of maximal dimension is
\begin{equation} P^u(t)(\Lambda(t)-\log|a|)^2 +
(P^u(t)+t\log|a|)\Lambda(t)^2=0.
\end{equation}
Therefore, the result follows from equation \RE{eqfordim}.\eproof
\begin{corollary}\label{corfinmesmax}
The set of all measures of maximal dimension is finite.
\end{corollary}
{\it Proof.} The function $\Delta$ is a non-constant real analytic
function on $[0,\infty)$. Therefore, it follows from the
uniqueness theorem for real analytic functions that $\Delta$ has
only finitely many maxima in $[t^s,t^u]$. The result follows from
Theorem \ref{mesmaxdim}.\eproof
\begin{corollary}
${\DD}(g)< \dim_H J$
\end{corollary}
{\it Proof.} Assume  $\DD(g)=\dim_H J$. By Theorem
\ref{mesmaxdim}, there exists a measure of maximal dimension.
Therefore, this measure is a measure of full dimension. But this
contradicts the standing assumption of this section that $g$ does
not  admit a measure of full dimension. \eproof
\begin{corollary}\label{corBer}
 Every measure $\nu$ of maximal  dimension is Bernoulli.
\end{corollary}
{\it Proof.} Since $g|_J$ is topological mixing (see \cite{BS1}),
the result follows from the fact that $\nu$ is an equilibrium
measure of a H\"older continuous potential (see \cite{B}, Thm.
4.1).\eproof
\section{Dependence on Parameters}
Let $A$ denote an open subset of $\bC^k$ and let
$(g_\lambda)_{\lambda\in A}$ be a holomorphic family of hyperbolic
mappings of the form $g_\lambda=g_{\lambda_1}\circ \cdots \circ
g_{\lambda_m}$, where each of the mappings $g_{\lambda_i}$ is a
generalized H\'enon mapping of fixed degree $d_i$ (see \cite{W1}
for more details). We denote by $J_\lambda$ the Julia set, by
$a_\lambda$ the Jacobian determinant,  and by $P^{u/s}_\lambda$
the unstable/stable pressure function of $g_\lambda$. We also
write $\Delta_\lambda(t)$ instead of $\Delta(t)$.  First, we show
 that $\DD(g)$ depends continuously on the
parameter of the mapping.
\begin{theorem}\label{thdepcon}
The function $\lambda\mapsto {\DD}(g_\lambda)$ is continuous in
$A$.
\end{theorem}
{\it Proof. } Let $\lambda_0\in A$. The result of \cite{VW}
implies that there exist $\epsilon
>0$ and a real analytic function
\begin{equation}
{\cal P}:B(\lambda_0,\epsilon)\times [0,\infty)\to \bR,
\end{equation}
such that ${\cal P}(\lambda,.)=P^u_\lambda$ for all $\lambda\in
B(\lambda_0,\epsilon)$. Therefore equations \RE{eq1},
\RE{eqfordim} imply that
\begin{equation}
{\cal D}:B(\lambda_0,\epsilon)\times [0,\infty)\to \bR,\quad
(\lambda,t)\mapsto \Delta_\lambda(t)
\end{equation}
is also a real analytic function. Now we may conclude by Corollary
\ref{thmefulun} and Theorem \ref{mesmaxdim} that
\begin{equation}
\DD(g_\lambda)= \max_{t\in [0,2]} {\cal D}(\lambda,t).
\end{equation}
The result follows  by an elementary continuity argument. \eproof
{ Remark. } McCluskey and Manning \cite{MM} considered $C^2$ Axiom
A diffeomorphisms of real surfaces. They showed that for these
mappings $\DD(g)$ depends continuously on the mapping with respect
to $C^2$ topology.\\[0.3cm]
 Finally, we show that  $\DD(g)$ depends
plurisubharmonically on the parameter of mapping.
\begin{theorem}\label{thdepsub}
The function $\lambda\mapsto {\DD}(g_\lambda)$ is plurisubharmonic
in $A$.
\end{theorem}
{\it Proof. } Let $g_0\in A$ and let $L$ be a complex line in
$\bC^k$ containing $g_0$. Then there exists a holomorphic family
$(g_\lambda)_{\lambda\in D}$, where $D$ is a disk with center $0$
in $\bC$ such that $\{g_\lambda:\lambda\in D\}$ is a neighborhood
of $g_0$ in $L\cap A$. If the radius of $D$ is small enough, then
there exists a family $(\kappa_\lambda)_{\lambda\in D}$, where
each $\kappa_\lambda$ is the uniquely defined conjugacy between
$g_0|_{J_0}$ and $g_\lambda|_{J_\lambda}$. Therefore,
$T_\lambda=(\kappa_\lambda)_*$ defines a family of bijections from
$M(J_0,g_0|_{J_0})$ to $M(J_\lambda,g_\lambda|_{J_\lambda})$.
Moreover we have $h_{\nu}(g_0)=h_{T_\lambda(\nu)}(g_\lambda)$ for
all $\nu\in M(J_0,g_0|_{J_0})$ and all $\lambda\in D$ (see
\cite{W1} for the details). In \cite{W1} we showed that  if
$\nu\in M(J_0,g_0|_{J_0})$ is fixed, then $\lambda\mapsto
\Lambda(T_\lambda(\nu))$ is a harmonic function in $D$.  We
conclude by Young's formula \RE{yo} that
\begin{equation}
\DD(g_\lambda)=\sup_{\nu\in M(J_0,g_0|_{J_0})}
\left(\frac{h_{\nu}(g_0)}{\Lambda(T_\lambda(\nu))}+
\frac{h_{\nu}(g_0)}{\Lambda(T_\lambda(\nu))-\log|a_\lambda|}
\right).
\end{equation}
The functions $\lambda\mapsto\Lambda(T_\lambda(\nu))$,
$\lambda\mapsto\Lambda(T_\lambda(\nu))-\log|a_\lambda|$ are
harmonic in $D$. Note that $x\mapsto x^{-1}$ is a convex function
on $\bR^+$. This implies that the functions $\lambda\mapsto
\Lambda(T_\lambda(\nu))^{-1}$, $\lambda\mapsto
(\Lambda(T_\lambda(\nu))-\log|a_\lambda|)^{-1}$ are subharmonic in
$D$. The continuous function $\lambda\mapsto \DD(g_\lambda)$ is
therefore given by the supremum over a family of subharmonic
functions. We conclude that the function $\lambda\mapsto
\DD(g_\lambda)$ is subharmonic in $D$. This completes the proof.
\eproof
\\[0.3cm] {\large \bf
Acknowledgement.} I would like to thank Eric Bedford and Marlies
Gerber for the helpful conversations during the preparation of
this paper.


\begin{thebibliography}{99}
{\small

\bibitem[B]{B}R. Bowen,  Equilibrium states and the ergodic
theory of Anosov diffeomorphisms, Lecture Notes in Math. {\bf
470}, Springer-Verlag, Berlin, 1975
\bibitem[BS1]{BS1}E. Bedford and  J. Smillie,  Polynomial diffeomorphisms of $\bC^2$ : Currents, equilibrium measure and hyperbolicity,  Invent.
Math. {\bf 103} (1991), 69 - 99
\bibitem[BS2]{BS2}E. Bedford and J. Smillie, Polynomial diffeomorphisms of $\bC^2$ 2: Stable manifolds and recurrence,  J. of the
AMS {\bf 4} (1991), 657 - 679
\bibitem[BS3]{BS3}E. Bedford and J. Smillie,  Polynomial diffeomorphisms of   $\bC^2$ 3:
 Ergodicity, exponents and entropy of the equilibrium measure,  Math.
 Ann.
 {\bf 294} (1992),
 395 - 420
\bibitem[BS6]{BS6}E. Bedford and J. Smillie,  Polynomial diffeomorphisms of   $\bC^2$ 6:
Connectivity of $J$,  Ann. of Math. {\bf 148} (1998), 695 - 735
\bibitem[DGS]{DGS}M. Denker, C. Grillenberger and K. Sigmund,
Ergodic theory on compact spaces,  Lecture Notes in Math. {\bf
527}, Springer-Verlag, Berlin, 1976
\bibitem[DU]{DU}M. Denker, M. Urbanski, On
Sullivan's conformal measures for rational maps on the Riemann
Sphere, Nonlinearity {\bf 4} (1991), 365-384
\bibitem[FM]{FM}S. Friedland and J. Milnor,  Dynamical properties of plane polynomial
automorphisms, Ergod. Theory and Dyn. Sys. {\bf 9} (1989),
 67-99
\bibitem[FO]{FO}S. Friedland and G. Ochs, Hausdorff dimension,
strong hyperbolicity and complex dynamics, Discrete and Continuous
Dyn. Syst. {\bf 3} (1998),
 405 - 430
   \bibitem[FS]{FS}J.E. Forn{\ae}ss and N. Sibony,  Complex H\'{e}non mappings in $\bC^2$ and Fatou-Bieberbach domains,
Duke Math. J.  {\bf 65} (1992), 345 - 380
\bibitem[GP]{GP}D. Gatzouras and Y. Peres, The variational
principle for Hausdorff dimension: a survey. Ergodic theory of
${\bf Z}^d$ actions (Warwick, 1993-1994), 113-125, London Math.
Soc. Lecture Note Ser., {\bf 228}, Cambridge Univ. Press,
Cambridge (1996)
 \bibitem[KH]{KH}A. Katok and B. Hasselblatt, Indroduction to
 modern theory of dynamical systems, Cambridge University
 Press, 1995
\bibitem[MM]{MM}H. McCluskey and A. Manning,  Hausdorff dimension
for Horseshoes, Ergod. Theory and Dyn. Syst. {\bf 3} (1983), 251 -
260
\bibitem[R]{R}D. Ruelle, Thermodynamic formalism, Addison-Wesley,
Reading, MA, 1978
\bibitem[U]{U}M. Urbanski, Measures and dimensions in conformal
dynamics, preprint, University of North Texas
\bibitem[VW]{VW}A. Verjovsky and H. Wu,  Hausdorff dimension of
Julia sets of complex H\'enon mappings,  Ergod. Theory and Dyn.
Syst {\bf 16} (1996), 849 - 861
\bibitem[Wa]{W}P. Walters,  An introduction to ergodic
theory, Springer, 1981
\bibitem[Wo1]{W1}C. Wolf,  Dimension of Julia sets of polynomial
automorphisms of $\bC^2$, Michigan Mathematical Journal {\bf 47}
(2000), 585 - 600
\bibitem[Wo2]{W2}C. Wolf,  Hausdorff and topological dimension for polynomial
automorphisms of $\bC^2$, to appear in Ergod. Theory and Dyn.
Syst.
\bibitem[Y]{Y}L.-S. Young, Dimension, entropy and Lyapunov
exponents,  Ergod. Theory and Dyn. Syst {\bf 2} (1982), 109-124 }
\end{thebibliography}
\end{document}